\documentclass[12pt]{article}
\usepackage{amscd, amssymb, latexsym, amsmath}

\newtheorem{thm}{Theorem}
\newtheorem{cor}{Corollary}
\newtheorem{pro}{Proposition}

\newtheorem{dfn}{Definition}

\newtheorem{rem}{Remark}

\newenvironment{proof}
{\noindent {\em Proof.}} {\hfill $\Box$}

\numberwithin{thm}{section} \numberwithin{cor}{section}
\numberwithin{pro}{section} \numberwithin{dfn}{section}
\numberwithin{lem}{section}
\numberwithin{rem}{section}\numberwithin{equation}{section}

\newcommand{\R}{\mathbb R}

\begin{document}
\title
{Isometric embeddings into the Minkowski space and new quasi-local mass}
\author{Mu-Tao Wang and Shing-Tung Yau}
\date{May 9, 2008, revised August 12, 2008}
\maketitle
\begin{abstract}
The definition of quasi-local mass for a bounded space-like region $\Omega$ in space-time is essential in several major unsettled problems in general relativity.
The quasi-local mass is expected to be a type of flux integral on the boundary two-surface $\Sigma=\partial \Omega$ and should be independent of whichever space-like region $\Sigma$ bounds. An important idea which is related to the Hamiltonian formulation of general relativity is to consider a reference surface in a flat ambient space with the same first fundamental form and derive the quasi-local mass from the difference of the extrinsic geometries. This approach has be taken by Brown-York \cite{by1}\cite{by2} and Liu-Yau \cite{ly}\cite{ly2} (see also related works \cite{hh}, \cite{la}, \cite{bly}, \cite{ki}, \cite{bm}, \cite{ep}, \cite{wy1}, \cite{za}) to define such notions using the isometric embedding theorem into the Euclidean three space. However, there exist surfaces in the Minkowski space whose quasilocal mass is strictly positive \cite{ost}. It appears that the momentum information needs to accounted for to reconcile the difference. In order to fully capture this information, we use isometric embeddings into the Minkowski space as references. In this article, we first prove an existence and uniqueness theorem for such isometric embeddings. We then solve the boundary value problem for Jang's \cite{ja} equation as a procedure to recognize such a surface in the Minkowski space. In doing so, we discover new expression of quasi-local mass for a large class of ``admissible" surfaces (see Theorem A and Remark 1.1). The new mass is positive when the ambient space-time satisfies the dominant energy condition and vanishes on surfaces in the Minkowski space. It also has the nice asymptotic behavior at spatial infinity and null infinity. Some of these results were announced in \cite{wy2}.
\end{abstract}

\section{Introduction}

\subsection{Dominant energy condition and positive mass theorem}
 Let $N$ be a space-time, i.e. a four-manifold  with a Lorentzian metric
$g_{\alpha\beta}$ of signature $(-+++)$ that satisfies the {\it
Einstein equation}:
\[R_{\alpha \beta}-\frac{s}{2}g_{\alpha \beta}=8\pi GT_{\alpha
\beta}\]
 where $R_{\alpha\beta}$ and $s$ are the Ricci curvature and the Ricci scalar curvature of $g_{\alpha\beta}$, respectively.
$G$ is the gravitational constant and $T_{\alpha\beta}$ is the energy-momentum
tensor of matter density. The metric $g_{\alpha\beta}$ defines space-like, time-like and null vectors on the tangent space of $N$ accordingly.

A ``dominant energy condition", which corresponds to a positivity condition on the matter density $T_{\alpha\beta}$, is expected to be satisfied on any realistic space-time. It means the following: for any time-like vector $e_0$, $T(e_{0}, e_0)\geq 0$ and $T(e_0,
\cdot)$ is a non-space-like co-vector. We shall assume throughout this article the space-time $N$ satisfies the dominant energy condition. Consider a space-like hypersurface $(M, g_{ij}, p_{ij})$ in $N$ where $g_{ij}$ is the induced (Riemannian) metric and
$p_{ij}$ is the second fundamental form with respect to the future-directed time-like unit normal vector field of $M$. The dominant energy condition together with the compatibility conditions for submanifolds imply
\begin{equation}\label{dec}\mu\geq |J|\end{equation} where \[\mu=\frac{1}{2}({R}-p_{ij}p^{ij}+(
p_{k}^k)^2),\] and
\[J^i= D_j(p^{ij}- p_{k}^{k} g^{ij}).\] Here
${R}$ is the scalar curvature of $M$.

An important special case is when $p_{ij}=0$ (time-symmetric case) and the dominant energy condition implies that the scalar curvature of $M$ is non-negative.

The positive mass theorem proved by Schoen-Yau \cite{sy1, sy2, sy3} (later a different proof by Witten \cite{wi} ) states:

\begin{thm}

Let $(M, g_{ij}, p_{ij})$ be a complete three manifold that satisfies (\ref{dec}). Suppose $M$ is {\it asymptotically flat}:
i.e. there exists a compact set $K\subset M$ such that $M\backslash K$ is diffeomorphic to a union of complements of balls in $\R^3$ (called ends) such that
 $g_{ij}=\delta_{ij} +a_{ij}$  with
$a_{ij}=O(\frac{1}{r})$, $\partial_k(a_{ij})=O(\frac{1}{r^2})$,
$\partial_l\partial_k(a_{ij})=O(\frac{1}{r^3}),$ and  $p_{ij}=O(\frac{1}{r^2})$,
$\partial_k(p_{ij})=O(\frac{1}{r^3})$ on each end of $M\backslash K$.

Then the ADM mass (Arnowitt-Deser-Misner) of each end of $M$  is positive, i.e.

\begin{equation}\label{pmt} E\geq |P|\end{equation} where

\[E=\lim_{r\rightarrow \infty}\frac{1}{16\pi G}
\int_{S_r}(\partial_j g_{ij}-\partial_i g_{jj})d\Omega^i\] is the total energy and

\[P_k=\lim_{r\rightarrow \infty} \frac{1}{16\pi G}\int_{S_r}
2(p_{ik}-\delta_{ik} p_{jj})d\Omega^i\] is the total momentum. Here $S_r$ is a coordinate sphere of radius $r$ on an end.  \end{thm}

We notice that the conclusion of the theorem is equivalent to the four-vector $(E, P_1, P_2,P_3)$ is
future-directed time-like, i.e.  \[E\geq 0 \,\,\,\text{and}\,\,
-E^2+P_1^2+P_2^2+P_3^2\leq 0.\]

The asymptotic flat condition can be considered a gauge condition to assure that $M$ can be compared to the flat space $\R^3$. The essence of the positive mass theorem is that positive local matter density (\ref{dec}) measured pointwise should imply positive total energy momentum (\ref{pmt}) measured at infinity. In contrast,  the ``quasi-local mass" corresponds to the measurement of mass of in-between scales.
\subsection{Two-surfaces in space-time and quasi-local notion of mass}
Let $N$ be a time-oriented space-time. Denote the Lorentzian metric on $N$ by $\langle \cdot, \cdot\rangle$ and covariant derivative by $\nabla^N$. Let $\Sigma$ be a closed space-like two-surface embedded in $N$. Denote the induced Riemannian metric on $\Sigma$ by $\sigma$ and the gradient and Laplace operator of $\sigma$ by $\nabla$ and $\Delta$, respectively.

 Given any two tangent vector $X$ and $Y$ of $\Sigma$, the second fundamental form of $\Sigma$ in $N$ is given by $\mbox{II}(X, Y)=(\nabla^N_{X} Y)^\perp$ where $(\cdot)^\perp$ denotes the projection onto the normal bundle of $\Sigma$. The mean curvature vector is the trace of the second fundamental form, or ${H}=tr_\Sigma \mbox{II}=\sum_{a=1}^2 \mbox{II}(e_a, e_a)$ where $\{e_1, e_2\}$ is an orthonormal basis of the tangent bundle of $\Sigma$.

 The normal bundle is of rank two
with structure group $SO(1,1)$ and the induced metric on the normal bundle is of signature $(-, +)$. Since the Lie algebra of $SO(1,1)$ is isomorphic to $\R$, the connection form of the normal bundle is a genuine 1-form that depends on the choice of the normal frames. The curvature of the normal bundle is then given by an exact 2-form which reflects the fact that any $SO(1,1)$ bundle is topologically trivial. Connections of different choices of normal frames differ by an exact form. We define
\begin{dfn}
Let $e_3$ be a space-like unit normal along $\Sigma$, the connection form determined by $e_3$ is defined to be
\begin{equation}\label{connection}\alpha_{e_3} (X)=\langle \nabla^N_X e_3, e_4\rangle \end{equation} where $e_4$  is the future-directed time-like unit normal that is orthogonal to $e_3$.

When $\Sigma$ bounds a space-like hypersurface $\Omega$ with $\partial \Omega=\Sigma$, we choose $e_3$ to be the space-like outward unit normal with respect to $\Omega$. The connection form is then denoted by $\alpha_\Omega$.
\end{dfn}

 Suppose $\Sigma$ bounds a space-like hypersurface $\Omega$ in $N$, the definition of quasi-local mass  $m_\Sigma$  asks that (see \cite{ea}, \cite{cy})

(1) $m_\Sigma\geq 0$ under the dominant energy condition.

(2) $m_\Sigma=0 $ if and only if $\Sigma$ is in the Minkowski spacetime.

(3) The limit of $m_\Sigma $ on large coordinates spheres of asymptotically flat (null) hypersurfaces should approach the ADM (Bondi) mass.

The quasi-local mass is supposed to be closely related to the formation of black holes according to the hoop conjecture of Throne. Various definitions for the quasi-local mass have been proposed (see for example the review article by Szabados \cite{sz}).


In this article, we shall focus on quasi-local mass defined by the following comparison principle: anchor the intrinsic geometry (the induced metric) by isometric embeddings and compare other extrinsic geometries. An important feature that we expect is the definition should be a flux type integral on $\Sigma$ and it should depend only on the fact that $\Sigma$ bounds a space-like hypersurface $\Omega$, but does not depend which specific $\Omega$ it bounds.

\subsection{Prior results}

We recall the solution of Weyl's isometric embedding problem by Nirenberg \cite{ni} and independently, Pogorelov \cite{po}:

\begin{thm}\label{isom}
Let $\Sigma$ be a closed surface with a Riemannian metric of positive Gauss curvature, then there exists an isometric embedding $i:\Sigma \hookrightarrow \R^3$ that is unique up to Euclidean rigid motions.
\end{thm}

In particular, the mean curvature of the isometric embedding is uniquely determined by the metric.  Through a Hamiltonian-Jacobi analysis of  Einstein's action, Brown and York \cite{by1} \cite{by2} introduced

\begin{dfn}
Suppose a two-surface $\Sigma$ bounds a space-like region $\Omega$ in a space-time $N$. Let $k_\Omega$ be the mean curvature
of $\Sigma$ with respect to the outward normal of $\Omega$. Assume the induced metric on $\Sigma$ has positive Gauss curvature and denote by $k_0$ the mean curvature of the isometric
embedding of $\Sigma$ into $\R^3$. The Brown-York mass  is defined to be:
\[\frac{1}{8\pi G}\left(\int_\Sigma k_0-\int_\Sigma k_\Omega\right).\]
\end{dfn}

Liu and Yau \cite{ly} \cite{ly2} (see also Kijowski \cite{ki}) defined

\begin{dfn}
Suppose $\Sigma$ is an embedded two-surface that bounds a space-like region in a space-time $N$. Assume $\Sigma$ has positive Gauss curvature. The Liu-Yau mass is defined to be
\[\frac{1}{8\pi G}\left(\int_\Sigma k_0-\int_\Sigma |H|\right)\] where $|H|$ is the Lorentzian norm of the mean curvature vector.
\end{dfn}
The Brown-York and Liu-Yau mass are proved to be positive by Shi-Tam \cite{st} in the time-symmetric case, and Liu-Yau \cite{ly} \cite{ly2}, respectively.
\begin{thm}\cite{st}
Suppose $\Omega$ has non-negative scalar curvature and $k_\Omega>0$. Then the Brown-York mass of $\Sigma$ is nonnegative and it equals zero if only if $\Omega$ is flat.
\end{thm}

\begin{thm}\cite{ly} \cite{ly2}
Suppose $N$ satisfies the dominant energy condition and the mean curvature vector of $\Sigma$ is space-like. The Liu-Yau mass is non-negative and it equals zero only if $N$ is isometric to $\R^{3,1}$ along $\Sigma$.
\end{thm}

However, \'{O} Murchadha, Szabados, and Tod \cite{ost} found examples of surfaces in the Minkowki space which satisfy the assumptions but whose Liu-Yau mass, as well as Brown-York, mass, are  strictly positive. It seems the missing of the momentum information $p_{ij}$ is responsible for this inconsistency: the Euclidean space can be considered as
 a totally geodesic space-like hypersurface in the Minkowski space with the second fundamental form $p_{ij}=0$ and in both the Brown-York and Liu-Yau case, the reference is taken to be the isometric embedding into $\R^3$. In order to capture the information of $p_{ij}$, we need to take the reference surface to be a general isometric embedding into the Minkowski space. However, an intrinsic difficulty for this embedding problem is that there are four unknowns (the coordinate functions in $\R^{3,1}$) but only three equations (for the first fundamental form). An ellipticity condition in replacement of the positive Gauss curvature condition is also needed to guarantee the uniqueness of the solution. We are able to achieve these in this article and indeed the extra unknown (corresponds to the time function) allows us to identify a canonical gauge in the physical space $N$ and define a quasi-local mass expression. We refer to our paper \cite{wy2} in which this expression was derived from the more physical point of view, i.e. the Hamilton-Jacobi analysis of the gravitational action.

\subsection{Results, organization and acknowledgement}
We first state the key comparison theorem:

\bigskip
\noindent\textbf{Theorem A}
\textit{
Let $N$ be a space-time that satisfies the dominant energy condition. Suppose $i:\Sigma \hookrightarrow N$ is a closed embedded space-like two-surface in $N$ with space-like mean curvature vector $H$. Let $i_0:\Sigma \hookrightarrow \R^{3,1}$ be an isometric embedding into the Minkowski space and let $\tau$ denote the restriction of the time function $t$ on $i_0(\Sigma)$. Let $\bar{e}_4$ be the future-directed time-like unit normal along $i(\Sigma)$ such that
\[\langle H, \bar{e}_4\rangle=\frac{-\Delta \tau}{\sqrt{1+|\nabla \tau|^2}}\] and $\bar{e}_3$ be the space-like unit normal along $\Sigma$ with $\langle \bar{e}_3, \bar{e}_4\rangle=0$ and $\langle H, \bar{e}_3\rangle <0$. Let $\widehat{\Sigma}$ be the projection of $i_0(\Sigma)$ onto $\R^3=\{t=0\}\subset \R^{3,1}$ and $\hat{k}$ be the mean curvature of $\widehat{\Sigma}$ in $\R^3$. If $\tau$ is admissible (see Definition \ref{adm}), then
\begin{equation}\label{pos_mass}\int_{\widehat{\Sigma}} \hat{k}-\int_\Sigma -\sqrt{1+|\nabla \tau|^2}\langle H, \bar{e}_3\rangle-\alpha_{\bar{e}_3} (\nabla \tau)\end{equation} is non-negative.}

\bigskip

Indeed, we show  \begin{equation}\label{total_mean}\int_{\widehat{\Sigma}} \hat{k}=\int_\Sigma -\sqrt{1+|\nabla \tau|^2}\langle H_0, \breve{e}_3\rangle-\alpha_{\breve{e}_3} (\nabla \tau)\end{equation} (see equation (\ref{int_hat_k}) ) where $H_0$ is the mean curvature vector of $i_0(\Sigma)$ in $\R^{3,1}$, $\breve{e}_3$ is the space-like unit normal along $i_0(\Sigma)$ in $\R^{3,1}$ that is orthogonal to the time direction. The expression (\ref{pos_mass}) naturally arises as the surface term in the Hamiltonian of gravitational action (see Remark \ref{surf_term}). When the reference isometric embedding lies in an $\R^3$ with $\tau=0$, it recovers the Liu-Yau mass.

\begin{rem}
If the Gauss curvature of $\Sigma$ is positive, an isometric embedding into an $\R^3$ with $\tau=0$ is admissible (see Corollary 5.3 and the preceding remark). In general, when the Gauss curvature is close to being positive, an isometric embedding with small enough $\tau$ is admissible.
\end{rem}

\begin{rem} We learned the expression in (\ref{total_mean}) from Gibbon's paper \cite{gi}. Indeed, we are motivated by \cite{gi} to study the projection of a space-like two-surface in the Minkowski space.
\end{rem}

The new quasi-local mass is defined to be the infimum of the expression  (\ref{pos_mass}) over all such isometric embeddings (see Definition \ref{ql_mass}). We prove that such embeddings are parametrized by the admissible $\tau$.

\bigskip
\noindent\textbf{Theorem B}
\textit{
Given a metric $\sigma$ and a function $\tau$ on $S^2$ such that the condition (\ref{convex_cond}) holds.  There
exists a unique space-like isometric embedding $i_0:S^2 \hookrightarrow \R^{3,1}$ with the induced metric $\sigma$ and the function $\tau$ as the time function.}
\bigskip

In \S 2, we study the expression $-\sqrt{1+|\nabla \tau|^2}\langle H, e_3\rangle-\alpha_{e_3} (\nabla \tau)$ for surfaces in space-time. We consider it as a generalized mean curvature and study the variation of the total integral. The gauge $\bar{e}_3, \bar{e}_4$ in Theorem A indeed minimizes the total integral (see Proposition 2.1). In \S 3, we prove Theorem B and study the total mean curvature of the projected surface. In particular, we prove equality (\ref{total_mean}). In \S 4, we study the boundary problem of Jang's equation and calculate the boundary terms. This is an important step in proving Theorem A. In \S 5, we define the new quasi-local mass and prove the positivity. In particular, Theorem A is proved. We emphasize that though the proof involves solving Jang's equation, the results depend only on the solvability but not on the specific solution. The Euler-Lagrange equation of the new quasi-local mass among all admissible $\tau$'s is derived in \S 6.
We wish to thank Richard Hamilton for helpful discussions on isometric embeddings and Melissa Liu for her interest and reading an earlier version of this  article. The first author would like to thank Naqing Xie for pointing out several typos in an earlier version.

\section{A generalization of mean curvature}
\begin{dfn}\label{little_j}
Suppose $i:\Sigma \hookrightarrow N$ is an embedded space-like two-surface. Given a smooth function $\tau$ on $\Sigma$ and a space-like normal $e_3$, the generalized mean curvature associated with these data is defined to be
\[h(\Sigma, i, \tau, e_3)=-\sqrt{1+|\nabla \tau|^2}\langle H, e_3\rangle-\alpha_{e_3} (\nabla \tau)\] where $H$ is the mean curvature vector of $\Sigma$ in $N$ and $\alpha_{e_3}$  is the connection form (see \ref{connection}) of the normal bundle of $\Sigma$ in $N$ determined by $e_3$ and the future-directed
 time-like unit normal $e_4$ orthogonal to $e_3$.
\end{dfn}

\begin{rem}\label{surf_term}
In the case when $\Sigma$ bounds a space-like region $\Omega$ and $e_3$ is the outward unit normal of $\Omega$, the mean curvature vector is
\[H=\langle H, e_3\rangle e_3-\langle H, e_4\rangle e_4.\]

We can reflect $H$ along the light cone of the normal bundle to get \[J=\langle H, e_4\rangle e_3-\langle H, e_3\rangle e_4.\]  Denote the tangent vector on $\Sigma$ dual to the one-form $\alpha_{e_3}$ by $V$, then the expression (3) in \cite{wy2} is $J-V$, where  $k=-\langle H, e_3\rangle$ and $p=-\langle H, e_4\rangle$. We have
\[h(\Sigma, i,\tau, e_3)=-\langle J-V, \sqrt{1+|\nabla \tau|^2} e_4-\nabla \tau\rangle .\] Notice that $\sqrt{1+|\nabla \tau|^2} e_4-\nabla \tau$ is again a future-directed unit time-like vector along $\Sigma$.
\end{rem}

Fix a base frame $\{\hat{e}_3, \hat{e}_4\}$ for the normal bundle, any other frame $\{e_3, e_4\}$ can be expressed as
\begin{equation}\label{normal_frames}e_3=\cosh\phi \hat{e}_3-\sinh\phi\hat{e}_4, \,\, e_4=-\sinh\phi \hat{e}_3+\cosh\phi \hat{e}_4\end{equation} for some $\phi$.

 We compute the integral \[\begin{split}&\int_{\Sigma} h(\Sigma, i, \tau, e_3) dv_\Sigma\\ &=\int_\Sigma \left[\sqrt{1+|\nabla \tau|^2}(\cosh \phi\langle\nabla^N_{e_a} \hat{e}_3, e_a\rangle-\sinh\phi \langle \nabla_{e_a}^N \hat{e}_4, e_a\rangle)-\alpha_{\hat{e}_3}(\nabla \tau)-\nabla \tau\cdot \nabla \phi \right] dv_\Sigma\end{split}\] and consider this expression as a functional of $\phi$.

 Suppose the mean curvature vector of $\Sigma$ is space-like, we may choose a base frame with
 \begin{equation}\label{hat_e3}\hat{e}_3=-\frac{H}{|H|}\end{equation} and $\hat{e}_4$ the future directed time-like unit normal that is orthogonal to $\hat{e}_3$.
 This choice makes $\langle \nabla_{e_a}^N \hat{e}_4, e_a\rangle=0$. Integration by parts, the functional becomes
\begin{equation}\int_\Sigma (\sqrt{1+|\nabla \tau|^2}\cosh \phi |H|-\alpha_{\hat{e}_3}(\nabla \tau)
+\phi \Delta \tau) dv_\Sigma.\end{equation}

As $|H|$ is positive, this is clearly a convex functional of $\phi$ which achieves the minimum as
\begin{equation}\label{sinh_phi}\sinh\phi=\frac{-\Delta \tau}{|H|\sqrt{1+|\nabla\tau|^2}}.\end{equation}
 We notice that the minimum is achieved by $\bar{e}_4$ such that
the expression \begin{equation}\label{bar_e4} |H|\sinh\phi=\langle H, \bar{e}_4\rangle= \frac{-\Delta \tau}{\sqrt{1+|\nabla \tau|^2}}\end{equation} depends only on $\tau$; this is taken as the characterizing property of $\bar{e}_4$ in \cite{wy2}.

\begin{dfn} \label{mathfrak_H} Given an isometric embedding $i:\Sigma\hookrightarrow N$ into a space-time with space-like mean curvature vector $H$. Denote
\[\mathfrak{H}(\Sigma, i, \tau)=\int_\Sigma \left[ \sqrt{(\Delta \tau)^2+|H|^2(1+|\nabla\tau|^2)}-\nabla \tau\cdot \nabla \phi-{\alpha}_{\hat{e}_3} (\nabla\tau)\right] dv_\Sigma.\]  where $\phi$ is defined by (\ref{sinh_phi}) and  ${\alpha}_{\hat{e}_3}$ is the connection one-form on $\Sigma$ associated with $\hat{e}_3$ in equation (\ref{hat_e3}). In terms of the frame $\bar{e}_3, \bar{e}_4$ where $\bar{e}_4$ is given by equation (\ref{bar_e4}) and $\bar{e}_3$ is the space-like unit normal with $\langle \bar{e}_3, \bar{e}_4\rangle=0$, then
\[\mathfrak{H}(\Sigma, i, \tau)=\int_\Sigma h(\Sigma, i, \tau, \bar{e}_3)dv_\Sigma=\int_\Sigma -\sqrt{1+|\nabla \tau|^2}\langle H, \bar{e}_3\rangle-\alpha_{\bar{e}_3} (\nabla \tau) dv_\Sigma. \]

\end{dfn}
\begin{pro}\label{can_gauge}
If the mean curvature vector of the embedding $i:\Sigma\hookrightarrow N$ is space-like and $e_3$ is any space-like unit normal such that $\langle H, e_3\rangle< 0$, then
\[\begin{split}& \int_\Sigma h(\Sigma, i, \tau, e_3) dv_\Sigma\geq  \mathfrak{H}(\Sigma, i, \tau). \end{split}\]
\end{pro}

\section{Isometric embeddings into the Minkowski space}

\subsection{Existence and uniqueness theorem}
Let $\Sigma$ be a two-surface diffeomorphic to $S^2$. We fix a Riemannian metric $\sigma$ on $\Sigma$, $\sigma=\sigma_{ab} du^a du^b$, in local coordinates $u^1, u^2$.  Denote the gradient, the Hessian, and the Laplace operator with respect to the metric $\sigma$ by $\nabla$, $\nabla^2$, and $\Delta$, respectively. We consider the isometric embedding problem of $(\Sigma, \sigma)$ into the Minkowski space $\R^{3,1}$ with prescribed mean curvature in a fixed time direction. Let $\langle\cdot, \cdot\rangle$ denote the standard Lorentzian metric on $\R^{3,1}$ and  $T_0$ be a constant unit time-like vector in $\R^{3,1}$, we have the following existence and uniqueness theorem:

\begin{thm}
Let  $\lambda$ be a function on $\Sigma$ with $\int_\Sigma \lambda dv_\Sigma=0$. Let $\tau$ be a potential function of $\lambda$, i.e. $\Delta \tau=\lambda$. Suppose
\begin{equation}\label{convex_cond}K+(1+|\nabla \tau|^2)^{-1}\det(\nabla^2 \tau)>0\end{equation} where $K$ is the Gauss curvature of $\sigma$ and $\det(\nabla^2\tau)$ is the determinant of the Hessian of $\tau$. Then there exists a unique space-like embedding
$X:\Sigma \hookrightarrow \R^{3,1}$ with the induced metric $\sigma$ and such that the mean curvature vector $H_0$ of the embedding satisfies \begin{equation}\label{prescribe_mean}\langle H_0, T_0\rangle =-\lambda.\end{equation}
\end{thm}

\begin{proof} We prove the uniqueness part first. Let $X_i:\Sigma \hookrightarrow \R^{3,1}$, $i=1,2$ be two isometric embeddings that satisfy (\ref{prescribe_mean}). Since the mean curvature vector of the embedding $X_i$ is $\Delta X_i$, this implies
\[\langle \Delta (X_1-X_2), T_0\rangle=0, \] or $\langle X_1-X_2, T_0 \rangle $ is a constant on $\Sigma$. Denote $\tau_i=-\langle X_i, T_0\rangle$, we thus have $d\tau_1=d\tau_2$. Now consider the projection $\widehat{X}_i:\Sigma \hookrightarrow \R^3$ onto the orthogonal complement of $T_0$; $\widehat{X}_i=X_i-\tau_iT_0$. The Gauss curvature of the embedding $\widehat{X}$ can be computed as
\begin{equation}\label{hat_K}\widehat{K}_i=(1+|\nabla \tau_i|^2)^{-1}[K+(1+|\nabla \tau_i|^2)^{-1}\det(\nabla^2 \tau_i)]\end{equation} which is positive by the assumption.

We compute the induced metric on the image of the embedding \[\langle d\widehat{X}_i, d\widehat{X}_i\rangle=\langle dX_i, dX_i\rangle +d\tau_i^2.\]

Since we assume $\langle dX_1, dX_1\rangle=\langle dX_2, dX_2\rangle =\sigma$, $\widehat{X}_i$'s are embeddings into $\R^3$ with the same induced metrics of positive Gauss curvature. By Theorem \ref{isom}, $\widehat{X}_1$ and $\widehat{X}_2$ are congruent in $\R^3$. Since $\tau_1$ and $\tau_2$ are different by a constant, $X_i$, as the graphs of $\tau_i$ over $\widehat{X}_i$, are congruent in $\R^{3,1}$.

We turn to the existence part. We start with the metric $\sigma$ and the function $\lambda$ and solve for $\tau$ in $\Delta \tau=\lambda$. The Gauss curvature $\widehat{K}$ of the new metric $\hat{\sigma}=\sigma+d\tau^2$ is again given by (\ref{hat_K}). Theorem \ref{isom} gives an embedding $\widehat{X}:\Sigma \hookrightarrow \R^3$ with the induced metric $\hat{\sigma}$. Now $X=\widehat{X}+\tau T_0 $ is the desired isometric embedding into $\R^{3,1}$ that satisfies (\ref{prescribe_mean}).
\end{proof}

The existence theorem can be formulated in terms of $\tau$ as the mean curvature vector is given by $H_0=\Delta X$.
\begin{cor} \label{time_coor} (\textbf{Theorem B})
Given a metric $\sigma$ and a function $\tau$ on $S^2$ such that the condition (\ref{convex_cond}) holds.  There
exists a unique space-like isometric embedding $i_0:S^2 \hookrightarrow \R^{3,1}$ with the induced metric $\sigma$ and the function $\tau$ as the time function.
\end{cor}
\subsection{Total mean curvature of the projection}

In this section, we compute the total mean curvature $\int_{\widehat{\Sigma}} \hat{k} dv_{\widehat{\Sigma}}$ of the projection $\widehat{\Sigma}$ in $\R^3$ in term of the geometry of $\Sigma$ in $\R^{3,1}$. Suppose $X:\Sigma\hookrightarrow \R^{3,1}$ is the embedding and $\tau=-\langle X, T_0\rangle$ is the restriction of the time function associated with $T_0$.
The outward unit normal $\hat{\nu}$ of $\widehat{\Sigma}$ in $\R^3$ and $T_0$ form an orthonormal basis for the normal bundle of $\widehat{\Sigma}$ in $\R^{3,1}$. Extend $\hat{\nu}$ along $T_0 $ by parallel translation and denote it by $\breve{e}_3$. We have
\begin{pro}
\begin{equation}\label{int_hat_k}\int_{\hat{\Sigma}}\hat{k} dv_{\widehat{\Sigma}}=\int_\Sigma \left[ -\langle H_0, \breve{e}_3\rangle \sqrt{1+|\nabla \tau|^2}-\alpha_{\breve{e}_3}(\nabla\tau) \right] dv_\Sigma\end{equation}
\end{pro}

\begin{proof} Denote by ${\nabla}^{\R^{3,1}}$ the flat connection associated with the Lorentzian metric on $\R^{3,1}$.
Take an orthonormal basis $\hat{e}_a, a=1, 2$ for the tangent space of $\widehat{\Sigma}$ and compute
\[\hat{k}=\langle {\nabla}^{\R^{3,1}}_{\hat{e}_a} \hat{\nu}, \hat{e}_a\rangle = \langle {\nabla}^{\R^{3,1}}_{\hat{e}_a} \hat{\nu}, \hat{e}_a\rangle+\langle {\nabla}^{\R^{3,1}}_{\hat{\nu}} \hat{\nu}, \hat{\nu}\rangle-\langle {\nabla}^{\R^{3,1}}_{T_0} \hat{\nu},T_0\rangle,\] because the last two terms are both zero.

Therefore $\hat{k}=g^{\alpha\beta} \langle{\nabla}^{\R^{3,1}}_{e_\alpha} \hat{\nu} , e_\beta\rangle$ for any orthonormal frame $e_\alpha$ of $\R^{3,1}$ where $g^{\alpha\beta}$ is the inverse of
$g_{\alpha\beta}= \langle e_\alpha, e_\beta\rangle$.

 Now $\breve{e}_3= \hat{\nu}$  may be considered as a space-like normal vector field along $\Sigma$. Pick an orthonormal basis $\{e_1, e_2\}$ tangent to $\Sigma$. Let $\breve{e}_4=\frac{1}{\sqrt{1+|\nabla\tau|^2}}(T_0-T_0^\top) $ be the future-directed unit normal vector in the direction of the normal part of $T_0$. It is not hard to see that $T_0^\top=-\nabla \tau$.  $\{\breve{e}_3, \breve{e}_4\}$ form an orthonormal basis for the normal bundle of $\Sigma$.
We derive
\begin{equation}\label{hat_k}\hat{k}=\langle{\nabla}^{\R^{3,1}}_{e_a}\breve{e}_3, e_a\rangle-\langle {\nabla}^{\R^{3,1}}_{\breve{e}_4}\breve{e}_3, \breve{e}_4\rangle=-\langle H_0, \breve{e}_3\rangle-\frac{1}{\sqrt{1+|\nabla \tau|^2}}\langle {\nabla}^{\R^{3,1}}_{\nabla \tau} \breve{e}_3, \breve{e}_4\rangle\end{equation} because $\hat{\nu}$ is extended along $T_0$ by parallel translation.

The area forms of $\Sigma$ and $\widehat{\Sigma}$ are related by  $dv_\Sigma=\frac{1}{\sqrt{1+|\nabla \tau|^2}}d v_{\widehat{\Sigma}}$. Integrating equation (\ref{hat_k}) over $\Sigma$, we obtain (\ref{int_hat_k})

\end{proof}

Suppose the mean curvature vector $H_0$ of $\Sigma$ in $\R^{3,1}$ is space-like. Let $e_3^{H_0}=\frac{-H_0}{|H_0|}$ be the unit vector in the direction of $H_0$ and $e_4^{H_0}$ the future-directed time-like unit normal vector with $\langle e_3^{H_0}, e_4^{H_0}\rangle=0$. Suppose that
\[e_3^{H_0}=\cosh\theta \breve{e}_3+\sinh\theta \breve{e}_4,\,\,\text{and}\,\,
e_4^{H_0}=\sinh\theta \breve{e}_3+\cosh\theta \breve{e}_4.\]

Since $\Delta\tau=-\langle H_0, T_0\rangle $ and $T_0=\sqrt{1+|\nabla\tau|^2}\breve{e}_4-\nabla \tau$, we derive  \begin{equation}\label{sinh_theta}\sinh\theta=\frac{-\Delta\tau}{|H_0|\sqrt{1+|\nabla\tau|^2}}.\end{equation} These imply the following relations \[\breve{e}_3 =\cosh\theta e_3^{H_0}-\sinh\theta e_4^{H_0}, \,\,\text{and}\,\,
\breve{e}_4=-\sinh\theta e_3^{H_0}+\cosh\theta e_4^{H_0}.\]

The integrand on the right hand side of (\ref{int_hat_k}) becomes
\[|H_0|\cosh\theta \sqrt{1+|\nabla\tau|^2}-\nabla\theta\cdot\nabla\tau-
\langle {\nabla}^{\R^{3,1}}_{\nabla\tau} e_3^{H_0},  e_4^{H_0}\rangle.\]

 Therefore we have
\begin{pro} When the mean curvature vector of $\Sigma$ in $\R^{3,1}$ is space-like, $\int_{\hat{\Sigma}}\hat{k} dv_{\widehat{\Sigma}}$ is equal to
\begin{equation}\label{int_hat_k2} \int_\Sigma \left[ \sqrt{(\Delta\tau)^2+|H_0|^2(1+|\nabla\tau|^2)}-\nabla \theta\cdot \nabla \tau-{\alpha}_{e_3^{H_0}} (\nabla\tau)\right] dv_\Sigma .\end{equation} where $\theta$ is given by (\ref{sinh_theta}) and
 $\alpha_{e_3^{H_0}}$ is the one-form on $\Sigma$ defined by $\alpha_{e_3^{H_0}}(X)=\langle {\nabla}^{\R^{3,1}}_{X} e_3^{H_0},  e_4^{H_0}\rangle$.
\end{pro}

\section{Jang's equation and boundary information}
\subsection{Jang's equation}
Jang's equation was proposed by Jang \cite{ja} in an attempt to solve the positive energy conjecture. Schoen and Yau came up with  different geometric interpretations, studied the equation in full, and applied to their proof \cite{sy3} of the positive mass theorem. Another important contribution of Schoen and Yau's work in \cite{sy3} is to understand the precise connection between the solvability of Jang's equation and the existence of black holes. This leads to the later works on the existence of black holes due to condensation of matter and boundary effect \cite{sy4} \cite{ya}. Given an initial data set $(\Omega, g_{ij}, p_{ij})$ where $p_{ij}$ is a symmetric two-tensor that represents the second fundamental form of $\Omega$ with respect to a future-directed time-like normal $e_4$ in a space-time $N$. We consider the Riemannian product $\Omega \times \R$ and extend $p_{ij}$ by parallel translation along the $\R$ direction to a symmetric tensor $P(\cdot, \cdot)$ on $\Omega\times \R$. Such an extension makes $P(\cdot, v)=0$ where $v$ denotes the downward unit vector in the $\R$ direction.

 Jang's equation asks for a hypersurface $\widetilde{\Omega}$ in $\Omega \times \R$, defined as  the graph of a function $f$ over $\Omega$, such that the mean curvature of $\widetilde{\Omega}$ in $\Omega\times \R$ is the same as the the trace of the restriction of $P$ to $\widetilde{\Omega}$. In terms of local coordinates $x^i$ on $\Omega$,
 the equation takes the form
\begin{equation}\label{jang1}\sum_{i, j=1}^3 (g^{ij}-\frac{f^if^j}{1+|D f|^2})(\frac{D_i D_j f}{(1+|D f|^2)^{1/2}}-p_{ij})=0, \end{equation}
where $Df=\frac{\partial f}{\partial x^i} g^{ij}\frac{\partial}{\partial x^j}$ is the gradient of $f$, $|Df|^2=g^{ij}\frac{\partial f}{\partial x^i}\frac{\partial f}{\partial x^j}$ and $D_iD_jf=\frac{\partial^2 f}{\partial x^i\partial x^j}
-\Gamma_{ij}^k\frac{\partial f}{\partial x^k}$ is the Hessian of $f$.

Pick an orthonormal basis $\{\tilde{e}_\alpha\}_{\alpha=1\cdots 4} $ for the tangent space of $\Omega\times \R$ along $\widetilde{\Omega}$ such that $\{\tilde{e}_i\}_{i=1\cdots3}$ is tangent to $\widetilde{\Omega}$ and $\tilde{e}_4$ is the downward unit normal, then Jang's equation is

\begin{equation}\label{jang2}\sum_{i=1}^3 \langle \widetilde{\nabla}_{\tilde{e}_i} \tilde{e}_4, \tilde{e}_i\rangle=\sum_{i=1}^3 P(\tilde{e}_i, \tilde{e}_i),\end{equation} here and throughout this section $\widetilde{\nabla}$ is the Levi-Civita connection on the product space $\Omega\times \R$.

\subsection{Boundary calculations}
  Let $\tau$ be a smooth function on $\Sigma=\partial \Omega$. We consider a solution $f$ of Jang's equation in $\Omega\times \R$ that satisfies the Dirichlet boundary condition $f=\tau$ on $\Sigma$.

  Denote the graph of $\tau$ over $\Sigma$ by $\widetilde{\Sigma}$ and the graph of $f$ over $\Omega$ by $\widetilde{\Omega}$ so that $\partial \widetilde{\Omega}=\widetilde{\Sigma}$. We choose orthonormal frames $\{e_1, e_2\}$ and $\{\tilde{e}_1, \tilde{e}_2\}$ for $T\Sigma$ and $T\widetilde{\Sigma}$, respectively. Let $e_3$ be the outward normal of $\Sigma$ that is tangent to $\Omega$. We also choose $\tilde{e}_3, \tilde{e}_4$ for the normal bundle of $\widetilde{\Sigma}$ in $\Omega\times \R$ such that $\tilde{e}_3$ is tangent to the graph $\widetilde{\Omega}$ and $\tilde{e}_4$ is a downward unit normal vector of $\widetilde{\Omega}$ in $\Omega\times \R$.  $\{e_1, e_2, e_3, v\}$ forms an orthonormal basis for the tangent space of $\Omega\times\R$, so does $\{\tilde{e}_\alpha\}_{\alpha=1\cdots 4}$. All these frames are extended along the $\R$ direction by parallel translation. Along $\Sigma$, we have
  \[D f=\nabla\tau+f_3 e_3\] where $f_3=e_3(f)$ is the normal derivative of $f$. $\tilde{e}_3$ and $\tilde{e}_4$ can be written down explicitly:

\begin{equation}\label{frames}\begin{split}\tilde{e}_3&=\frac{1}{\sqrt{1+|D f|^2}}\left[\sqrt{1+|\nabla \tau|^2} e_3-\frac{f_3}{\sqrt{1+|\nabla \tau|^2}}(v+\nabla \tau)\right]\,\,\text{and}\\
\tilde{e}_4&=\frac{1}{\sqrt{1+|D f|^2}}( v+D f).\end{split}
\end{equation} We check that $\tilde{e}_3$ and $\tilde{e}_4$ are orthogonal to $e_a-e_a(\tau)v$ for $a=1,2$.

Simple calculations yield
\begin{equation}\label{inner_products}\begin{split}
\langle e_3, \tilde{e}_3\rangle =\frac{\sqrt{1+|\nabla \tau|^2}}{\sqrt{1+|D f|^2}}, \,\, \text{and}\,\,
\langle e_3, \tilde{e}_4\rangle=\frac{f_3}{\sqrt{1+|D f|^2}}.\end{split}\end{equation}

Let $\tilde{k}=\langle \widetilde{\nabla}_{\tilde{e}_a} \tilde{e}_3,\tilde{e}_a\rangle$ be the mean curvature of $\widetilde{\Sigma}$ with respect to $\widetilde{\Omega}$. We are particularly interested in the following expression on $\widetilde{\Sigma}$:

\begin{equation}\label{til_k} \tilde{k}-\langle\widetilde{\nabla}_{\tilde{e}_4} \tilde{e}_4, \tilde{e}_3\rangle+P(\tilde{e}_4, \tilde{e}_3).\end{equation}

\begin{thm}\label{jang_boundary}
Let $i:\Sigma \hookrightarrow N $ be a space-like  embedding.
 Given any smooth function $\tau$ on $\Sigma$ and any space-like hypersurface $\Omega$ with $\partial\Omega=\Sigma$. Suppose the Dirichlet problem of Jang's equation (\ref{jang1}) over $\Omega$ subject to the boundary condition that $f=\tau$ on $\Sigma$ is solvable.  Then there exists a space-like unit normal $e'_3$ along $\Sigma$ in $N$ such that the expression (\ref{til_k})  at $\tilde{q}\in \widetilde{\Sigma}$ is equal to \[-\langle H, e'_3\rangle-(1+|\nabla \tau|^2)^{-1/2} \alpha_{e'_3}(\nabla \tau)\text{  at  } q\in \Sigma,\] where $\tilde{q}=(q, \tau(q))\in \widetilde{\Sigma}.$ In particular
\begin{equation}\label{int_tilde_k} \int_{\widetilde{\Sigma}} \tilde{k}-\langle\widetilde{\nabla}_{\tilde{e}_4} \tilde{e}_4, \tilde{e}_3\rangle+P(\tilde{e}_4, \tilde{e}_3) dv_{\widetilde{\Sigma}}=\int_\Sigma -\sqrt{1+|\nabla \tau|^2} \langle H, e'_3\rangle- \alpha_{e'_3}(\nabla \tau) dv_\Sigma.\end{equation}

Let $e_3$ be the outward unit normal of $\Sigma$ that is tangent to $\Omega$ and $e_4$ is the future-directed time-like normal of $\Omega$ in $N$, $e_3'$ is given by
\begin{equation}\label{cosh}e'_3=\cosh \phi e_3+\sinh \phi e_4 , \,\,\,\text{where}\,\,\,  \sinh\phi=\frac{-
f_3}{\sqrt{1+|\nabla \tau|^2}}.\end{equation}
\end{thm}

\begin{proof} The proof is through a sequence of calculations using the product structure of $\Omega\times \R$ and Jang's equation. It also relies on the fact that $P$, $\{\tilde{e}_\alpha\}_{\alpha=1}^4$, and $\{e_1, e_2, e_3, v\}$ are all parallel in the direction of $v$.  We first prove the following identity:
\begin{equation}\begin{split}\label{tilde_k}&\tilde{k}-\langle\widetilde{\nabla}_{\tilde{e}_4} \tilde{e}_4, \tilde{e}_3\rangle+P(\tilde{e}_3, \tilde{e}_4)\\
&=\langle \widetilde{\nabla}_{e_a}\tilde{e}_3, e_a\rangle+\frac{\langle e_3, \tilde{e}_4\rangle}{\langle e_3, \tilde{e}_3\rangle}\langle \widetilde{\nabla}_{e_a}\tilde{e}_4, e_a\rangle
-\frac{\langle e_3, \tilde{e}_4\rangle}{\langle e_3, \tilde{e}_3\rangle} P(e_a, e_a)+\frac{1}{\langle e_3, \tilde{e}_3\rangle}P(e_3, \tilde{e}_4-\langle \tilde{e}_4, e_3\rangle e_3) .\end{split}\end{equation}

We compute the terms $\langle\widetilde{\nabla}_{e_a} \tilde{e}_3, e_a\rangle$ and $\langle\widetilde{\nabla}_{e_a} \tilde{e}_4, e_a\rangle$ in the following. \[\langle\widetilde{\nabla}_{e_a}\tilde{e}_3, e_a\rangle=\sum_{i=1}^3 \langle\widetilde{\nabla}_{e_i} \tilde{e}_3, e_i \rangle-\langle\widetilde{\nabla}_{e_3} \tilde{e}_3, e_3\rangle=\sum_{\alpha=1}^4\langle\widetilde{\nabla}_{\tilde{e}_{\alpha}} \tilde{e}_3, \tilde{e}_\alpha\rangle-\langle\widetilde{\nabla}_{e_3} \tilde{e}_3, e_3\rangle,\] as $\{e_1, e_2, e_3, v\}$ and $\{\tilde{e}_{\alpha}\}_{\alpha=1\cdots 4}$ are both orthonormal frames for the tangent space of $\Omega\times \R$ and $\widetilde{\nabla}_v \tilde{e}_3=0$.

Notice that  $\langle\widetilde{\nabla}_{\tilde{e}_{3}} \tilde{e}_3, \tilde{e}_3\rangle=0$ and thus we obtain

\begin{equation}\label{tilde_e_3}
\sum_{a=1}^2\langle \widetilde{\nabla}_{e_a}\tilde{e}_3, e_a\rangle =\tilde{k}+\langle \widetilde{\nabla}_{\tilde{e}_4} \tilde{e}_3, \tilde{e}_4\rangle-\langle \widetilde{\nabla}_{e_3} \tilde{e}_3, e_3\rangle.\end{equation}

On the other hand,
\[\sum_{a=1}^2\langle\widetilde{\nabla}_{e_a}\tilde{e}_4, e_a\rangle=\sum_{i=1}^3\langle\widetilde{\nabla}_{e_i} \tilde{e}_4, e_i\rangle-\langle\widetilde{\nabla}_{e_3} \tilde{e}_4, e_3\rangle
=\sum_{i=1}^3\langle \widetilde{\nabla}_{\tilde{e}_i}\tilde{e}_4, \tilde{e}_i\rangle-\langle\widetilde{\nabla}_{e_3} \tilde{e}_4, e_3\rangle.\]

Applying Jang's equation (\ref{jang2}), we obtain
\[\sum_{a=1}^2\langle\widetilde{\nabla}_{e_a}\tilde{e}_4, e_a\rangle=\sum_{i=1}^3P(\tilde{e}_i, \tilde{e}_i)-\langle\widetilde{\nabla}_{e_3} \tilde{e}_4, e_3\rangle.\]

Furthermore, we derive
\[\sum_{i=1}^3P(\tilde{e}_i, \tilde{e}_i)=\sum_{i=1}^3 P(e_i, e_i)
+\frac{\langle e_3, \tilde{e}_3\rangle}{\langle e_3, \tilde{e}_4\rangle} P(\tilde{e}_3, \tilde{e}_4)-\frac{1}{\langle e_3, \tilde{e}_4\rangle} P({e}_3, \tilde{e}_4),\] using
\[\sum_{i=1}^3P(\tilde{e}_i, \tilde{e}_i)=\sum_{\alpha=1}^4 P(\tilde{e}_\alpha, \tilde{e}_\alpha)-P(\tilde{e}_4, \tilde{e}_4)=\sum_{i=1}^3P({e}_i, {e}_i)-P(\tilde{e}_4, \tilde{e}_4)\] and \[\langle e_3, \tilde{e}_4\rangle P(\tilde{e}_4, \tilde{e}_4)=P(e_3-\langle e_3,\tilde{e}_3\rangle\tilde{e}_3 ,\tilde{e}_4).\]

Therefore, we arrive at
\begin{equation}\label{tilde_e_4}\begin{split}\sum_{a=1}^2\langle\widetilde{\nabla}_{e_a}\tilde{e}_4, e_a\rangle &=\sum_{a=1}^2P({e}_a,{e}_a)+\frac{\langle e_3, \tilde{e}_3\rangle}{\langle e_3, \tilde{e}_4\rangle} P(\tilde{e}_3, \tilde{e}_4)\\
&-\frac{1}{\langle e_3, \tilde{e}_4\rangle} P({e}_3, \tilde{e}_4-\langle \tilde{e}_4, e_3\rangle e_3)-\langle \widetilde{\nabla}_{e_3}\tilde{e}_4, e_3\rangle .\end{split}\end{equation}
Combining (\ref{tilde_e_3}) and (\ref{tilde_e_4}) yields (\ref{tilde_k}).

Let $k=\langle \widetilde{\nabla}_{e_a} e_3, e_a\rangle $ be the mean curvature of $\Sigma$ (as the boundary of $\Omega$) with respect to $e_3$. As $\langle e_3, \tilde{e}_a\rangle=0$ for $a=1, 2$, we have  $e_3=\langle e_3, \tilde{e}_3\rangle \tilde{e}_3+\langle e_3, \tilde{e}_4\rangle \tilde{e}_4$. Plug this into the expression for $k$, we obtain
\[\begin{split}k&=\langle e_3, \tilde{e}_3\rangle \langle\widetilde{\nabla}_{e_a} \tilde{e}_3, e_a\rangle+\langle e_3, \tilde{e}_4\rangle \langle\widetilde{\nabla}_{e_a} \tilde{e}_4, e_a\rangle\\
&-\langle e_3, \tilde{e}_3\rangle e_a(\langle e_a, \tilde{e}_3\rangle)-\langle e_3, \tilde{e}_4\rangle e_a(\langle e_a,\tilde{e}_4\rangle).\end{split}\]

From here we solve for $\langle \widetilde{\nabla}_{e_a}\tilde{e}_3, e_a\rangle+\frac{\langle e_3, \tilde{e}_4\rangle}{\langle e_3, \tilde{e}_3\rangle}\langle \widetilde{\nabla}_{e_a}\tilde{e}_4, e_a\rangle$ and substitute into (\ref{tilde_k}) to obtain
\begin{equation}\label{tilde_k2}\begin{split}&\tilde{k}-\langle\widetilde{\nabla}_{\tilde{e}_4} \tilde{e}_4, \tilde{e}_3\rangle+P(\tilde{e}_3, \tilde{e}_4)\\
&=\frac{1}{\langle e_3, \tilde{e}_3\rangle}k-\frac{\langle e_3, \tilde{e}_4\rangle}{\langle e_3, \tilde{e}_3\rangle} P(e_a, e_a)+\frac{1}{\langle e_3, \tilde{e}_3\rangle}P(e_3, \tilde{e}_4-\langle \tilde{e}_4, e_3\rangle e_3) \\
&+e_a(\langle e_a, \tilde{e}_3\rangle)+\frac{\langle e_3, \tilde{e}_4\rangle}{\langle e_3, \tilde{e}_3\rangle} e_a(\langle e_a,\tilde{e}_4\rangle).\end{split}\end{equation}

We calculate the right hand side of (\ref{tilde_k2}) using (\ref{frames}) and

\[P(e_3, \tilde{e}_4-\langle \tilde{e}_4, e_3\rangle e_3)=\frac{1}{\sqrt{1+|D f|^2}}P(e_3, \nabla \tau).\]

The last two terms can also be calculated using (\ref{frames}) and
\[\begin{split} &e_a(\langle e_a, \tilde{e}_3\rangle)+ \frac{\langle e_3, \tilde{e}_4\rangle}{\langle e_3, \tilde{e}_3\rangle}e_a(\langle e_a,\tilde{e}_4\rangle)\\
&=div_\Sigma (\frac{-f_3 \nabla \tau}{\sqrt{(1+|D f|^2)(1+|\nabla \tau|^2)}})+\frac{f_3}{\sqrt{1+|\nabla\tau|^2}}div_\Sigma(\frac{\nabla \tau}{\sqrt{1+|D f|^2}})\\
&=-\frac{1}{\sqrt{1+|D f|^2}}\nabla \tau \cdot \nabla (\frac{f_3}{\sqrt{1+|\nabla \tau|^2}}).\end{split}\]

Recalling the definition of $\phi$ from (\ref{cosh}), this is equal to \[\frac{\nabla \tau\cdot \nabla \phi}{\sqrt{1+|\nabla \tau|^2}}.\]

The right hand side of (\ref{tilde_k2}) is therefore
\begin{equation}\label{momentum}\begin{split}
(1+|\nabla\tau|^2)^{-1/2}\left[\sqrt{1+|D f|^2}k-f_3 P(e_a, e_a)+P(e_3, \nabla \tau)
+\nabla \tau\cdot \nabla \phi\right] \end{split}\end{equation}

This is an expression on $\Sigma$ that depends on the functions $\tau$ and $f_3$ on $\Sigma$. Recall that the symmetric tensor $P$ originates from the second fundamental form of $\Omega$ with respect to the future-directed unit time-like normal $e_4$ in the space-time $N$.  Rewrite the expression (\ref{momentum}) in terms of $e_3$ and $e_4$:

\begin{equation}\label{tilde_k3}\begin{split}&\tilde{k}-\langle\widetilde{\nabla}_{\tilde{e}_4} \tilde{e}_4, \tilde{e}_3\rangle+P(\tilde{e}_3, \tilde{e}_4)\\
&=(1+|\nabla\tau|^2)^{-1/2}\left[\sqrt{1+|D f|^2}\langle \nabla^N_{e_a} e_3, e_a\rangle-f_3 \langle \nabla^N_{e_a} e_4, e_a\rangle -\alpha_{e_3}(\nabla \tau)+\nabla \tau\cdot \nabla \phi\right]
\end{split}\end{equation}

 On the other hand, with the orthonormal frame $e'_3, e'_4$ given by
 \[\label{framerelation}e'_3=\cosh \phi e_3+\sinh \phi e_4,\,\,e'_4=\sinh \phi e_3+\cosh \phi e_4 ,\] we compute
\[\begin{split}&\langle \nabla^N_{e_a} e'_3, e_a\rangle-(1+|\nabla \tau|^2)^{-1/2}\alpha_{e'_3} (\nabla \tau)\\
&=\cosh\phi \langle \nabla^N_{e_a} e_3, e_a\rangle+\sinh\phi \langle \nabla^N_{e_a} e_4, e_a\rangle -(1+|\nabla \tau|^2)^{-1/2}(\alpha_{e_3}(\nabla \tau)-\nabla \tau\cdot \nabla \phi).\end{split}\]

Plug in the expression for $\cosh\phi$ and $\sinh\phi$, we recover the right hand side of (\ref{tilde_k3}).

\end{proof}

\subsection{Boundary gradient estimate}
In this section, we demonstrate a sufficient condition for Jang's equation to be solvable. As most estimates are derived in Schoen-Yau's original paper \cite{sy3} for the asymptotically flat case, it suffices to control the boundary gradient of the solution.

\begin{thm}
The normal derivative of a solution of the Dirichlet problem of Jang's equation is bounded if $k >|tr_{\widetilde{\Sigma}} P|$.
\end{thm}

\begin{proof}
We consider the operator

\[Q(f)=(g^{ij}-\frac{f^if^j}{1+|D f|^2})\frac{D_i D_j f}{(1+|D f|^2)^{1/2}}-tr_{\widetilde{\Omega}} P .\]
where $\widetilde{\Omega}$ is the graph of $f$ over $\Omega$. The point is to construct sub and super solutions of this operator with the prescribed boundary condition.

Denote by $d$ the distance function to $\partial \Omega$. We extend the boundary data $\tau$ to the interior of $\Omega$, still denoted by $\tau$.  Consider a test function $f=\psi (d)+ \tau$ as the one in (14.11) of \cite{gt} where $\psi (d)=\frac{1}{\nu} \log (1+\kappa d)$ with $\kappa,\nu>0$. In particular $\psi''=-\nu(\psi')^2<0$, $\psi>0$, and $\psi'(d) \rightarrow \infty$ as $\kappa\rightarrow \infty$.  We compute
\[D_i D_j f=\psi'' d_i d_j+\psi' D_i D_j d+D_i D_j \tau.\]

Therefore,
\[\begin{split}&(g^{ij}-\frac{f^if^j}{1+|D f|^2})\frac{D_i D_j f}{(1+|D f|^2)^{1/2}}\\
&=\psi''(g^{ij}-\frac{f^if^j}{1+|D f|^2})\frac{d_i d_j}{(1+|D f|^2)^{1/2}}+\psi'(g^{ij}-\frac{f^if^j}{1+|D f|^2})\frac{D_i D_j d}{(1+|D f|^2)^{1/2}}\\
&+(g^{ij}-\frac{f^if^j}{1+|D f|^2})\frac{D_i D_j \tau}{(1+|D f|^2)^{1/2}}\end{split}.\]

Applying
\[\frac{1}{1+|D f|^2} g^{ij}\leq g^{ij}-\frac{f^i f^j}{1+|D f|^2}\leq g^{ij}\]
and
\[ |D d|=1, \] We derive the first term is bounded above by
\[\frac{\psi''}{(1+|D f|^2)^{3/2}}\] and the third term is bounded above by
\[\frac{|D^2\tau|}{(1+|D f|^2)^{1/2}}.\]

The second term is
\[\psi'(g^{ij}-\frac{f^if^j}{1+|D f|^2})\frac{D_i D_j d}{(1+|D f|^2)^{1/2}}=\psi'\frac{\Delta d}{(1+|D f|^2)^{1/2}}-\frac{\psi'}{(1+|D f|^2)^{3/2}}f^i f^j D_i D_j d.\]

We compute
\[ f^i f^j D_i D_j d=(\psi' d^i+\tau^i)(\psi' d^j+\tau^j) D_i D_j d=\tau^i\tau^j D_i D_j d\] where we used the identity
\[d^i D_i D_j d=0.\]

Therefore, $Q(f)$ is bounded from above by
\[\frac{\psi''+(1+|D f|^2)|D^2\tau|-\psi'\tau^i\tau^j D_i D_j d}{(1+|D f|^2)^{3/2}}+\psi'\frac{\psi' \Delta d}{(1+|D f|^2)^{1/2}}-tr_{\widetilde{\Omega}} P.\]

We recall that $\widetilde{\Omega}$ is the graph of $\psi(d)+\tau$ over $\Omega$. Let $\Omega_a=\{d\leq a\}\cap \Omega$ and $\widetilde{\partial \Omega_a}$ be the graph of $\psi(d)+\tau$ over $\partial \Omega_a$. We have
\[tr_{\tilde{\Omega}}P=tr_{\widetilde{\partial \Omega_a}} P+\frac{P(D d, D d)}{\sqrt{1+(\psi'+D \tau\cdot D d)^2}}\]

 Therefore $Q(f)$ is bounded from above by
\[\begin{split}
 &\frac{\psi''+(1+|D f|^2)|D^2\tau|-\psi'\tau^i\tau^j D_i D_j d}{(1+|D f|^2)^{3/2}}+\frac{\psi' \Delta d}{(1+|D f|^2)^{1/2}}-tr_{\widetilde{\partial \Omega_a}} P-\frac{P(D d, D d)}{\sqrt{1+(\psi'+D \tau\cdot D d)^2}}\end{split}.\]

 When $\tau=0$, this recovers formula (5.11) in \cite{ya}. In general, we recall $D f=\psi' D d+D \tau$ and

 \[|D f|^2\geq \theta (\psi')^2-\frac{\theta}{1-\theta}|D \tau|^2\] for any positive $\theta<1$.
We notice that $\Delta d$ approaches $-k$, where $k$ is the mean curvature of $\partial\Omega$ in $\Omega$. However, $tr_{\widetilde{\partial \Omega_a}} P$ approaches $tr_{\widetilde{\Sigma}} P$. Thus a sub and a super solution exist when $k\geq |tr_{\widetilde{\Sigma}} P|$.

 \end{proof}

\section{New quasi-local mass and the positivity}

First we define an admissible time function for a surface in space-time.

\begin{dfn}\label{adm}
Given a space-like embedding $i:\Sigma\hookrightarrow N$. A smooth function $\tau$ on $\Sigma$ is said to be admissible if

\noindent (1)\,\,$K+(1+|\nabla \tau|^2)^{-1}\det(\nabla^2\tau)>0$.

\noindent (2)\,\,$\Sigma$ bounds an embedded space-like three-manifold $\Omega$ in $N$ such that Jang's equation (\ref{jang1}) with the Dirichlet boundary data $\tau$ is solvable on $\Omega$.

\noindent (3)\,\,The generalized mean curvature $h(\Sigma, i, \tau, e_3')>0$ for the space-like unit normal $e_3'$ (\ref{cosh}) determined by Jang's equation.
\end{dfn}
We are now ready to define the quasi-local mass.

\begin{dfn}\label{ql_mass}
Given a space-like embedding $i:\Sigma\hookrightarrow N$. Suppose the set of admissible functions is non-empty. The quasi-local mass is the defined to be the infimum of
\[ \mathfrak{H}(\Sigma, i_0, \tau)-\mathfrak{H}(\Sigma, i, \tau)\] among all admissible $\tau$ where $\mathfrak{H}$ is given by Definition \ref{mathfrak_H} and $i_0$ is the unique space-like isometric embedding of $\Sigma$ into $\R^{3,1}$ associated with $\tau$ given by Theorem B.
\end{dfn}

The proof of the positivity of quasi-local mass is based on the following theorem which can be considered as a total mean curvature comparison theorem for solutions of Jang's equation.

\begin{thm}\label{jang_sol}
Suppose $\Omega$ is a Riemannian three-manifold with boundary $\Sigma$ and suppose there exists a vector field $X$ on $\Omega$ such that \begin{equation}\label{jang_scalar} R\geq 2|X|^2-2div X\end{equation} in $\Omega$ where $R$ is the scalar curvature of $\Omega$
and \begin{equation}\label{jang_mean} k>\langle X, \nu\rangle\end{equation} on $\Sigma$, where $\nu$ is outward normal of $\Sigma$ and $k$ is the
mean curvature of $\Sigma$ with respect to $\nu$. Suppose the Gauss curvature of $\Sigma$ is positive and $k_0$ is the mean curvature of the isometric embedding of $\Sigma$ into $\R^3$. Then

\[\int_\Sigma k_0 dv_\Sigma \geq \int_\Sigma k-\langle X, \nu\rangle dv_\Sigma.\]
\end{thm}

\begin{rem}\label{shitam}When $X=0$, the theorem was proved by Shi-Tam \cite{st}. By the calculation in Schoen-Yau \cite{sy2}, the condition  (\ref{jang_scalar}) holds for any solution of Jang's equation over an initial data set that satisfies the dominant energy condition.  The vector field $X$ is the dual of $\langle \widetilde{\nabla}_{\tilde{e}_4} \tilde{e}_4, \cdot\rangle-P(\tilde{e}_4, \cdot)$ in the notation of \S 4. In this case, Liu-Yau \cite{ly} essentially proved the theorem by conformally changing the metric to zero scalar curvature. The proof of Theorem 6.2 in \cite{wy1} gives a direct proof without conformal change in a slightly different setting.
\end{rem}

\begin{proof} The idea of the proof is similar to the one by Shi-Tam. Consider the isometric embedding of $\Sigma$ into $\R^3$ and denote the region inside the image $\Sigma_0$ by $\Omega_0$. We then glue together $\Omega$ and $\R^3\backslash \Omega_0$ along the identification of $\Sigma$ and $\Sigma_0$. Write the metric on $\R^3\backslash \Omega_0$ into the form $dr^2+ g_r$ where $r$ is the distance function to $\Sigma_0$ and $g_r$ is the induced metric on the level set $\Sigma_r$ of $r$. Applying Bartnik's \cite{ba2} quasi-spherical construction, we consider a new metric on $\R^3\backslash \Omega_0$ of the form $u^2 dr^2+g_r$ with zero scalar curvature and $u=\frac{k_0}{k-\langle X, \nu\rangle}$ at $r=0$. $u$ then satisfies a parabolic equation and the solution gives an asymptotically flat metric on $\R^3\backslash \Omega_0$. Denote by $\tilde{M}$ the space $\Omega\cup_\Sigma \R^3\backslash \Omega_0$ with the new metric $u^2dr^2+g_r$ on $\R^3\backslash \Omega_0$. The initial condition on $u$ implies the mean curvature of $\Sigma_0$ with respect to this new metric is $k-\langle X, \nu\rangle$. We still have the monotonicity formula, i.e. $\frac{d}{dr}\int_{\Sigma_r} k_0(r)(1-\frac{1}{u}) dv_{\Sigma_r} \leq 0$ where $k_0(r)$ is the mean curvature of $\Sigma_r$ in $\R^3$. Therefore, it remains to prove the positivity of the total mass of $\tilde{M}$. In the following, we prove a Lichnerowicz formula for such a manifold and  the existence of  harmonic spinors asymptotic to constant spinors.  According to the standard Lichnerowicz formula, on $\Omega$ we have
\begin{equation}\label{killing}\begin{split}&\int_\Omega
|\nabla\psi|^2+\frac{1}{4}\int_\Omega
R |\psi|^2-\int_\Omega |D\psi|^2\\
&=\int_{\partial \Omega} \langle \psi,
{\nabla}_{\nu}\psi +c(\nu) \cdot
{D}\psi \rangle,
\end{split}
\end{equation} where $\psi$ is a spinor, $c(\cdot)$ is the Clifford multiplication, $\nabla$ is the spin connection, and $D$ is the Dirac operator.

Integrating by parts, we obtain
\[\int_{\partial \Omega}\langle X, \nu\rangle
|\psi|^2=\int_\Omega div X
|\psi|^2+\int_\Omega X(|\psi|^2).\]

Formula (\ref{killing}) is equivalent to
\begin{equation}\begin{split}&\int_\Omega
|\nabla\psi|^2+\frac{1}{4}\int_\Omega (R+2div
X)|\psi|^2+\frac{1}{2}\int_\Omega X(|\psi|^2)-\int_\Omega |D\psi|^2\\
&=\int_{\partial \Omega} \langle \psi,
{\nabla}_{\nu}\psi+c(\nu) \cdot
{D}\psi \rangle+\frac{1}{2} \int_{\partial \Omega}\langle
X, \nu\rangle |\psi|^2.
\end{split}
\end{equation}

The boundary term can be rearranged as
\[\int_{\partial
\Omega}-\langle \psi, D^{\partial
\Omega}\psi\rangle-\frac{1}{2}(k-\langle X,
\nu\rangle)|\psi|^2,\] where $-D^{\partial
\Omega}\psi=\sum_{a=1}^2 c(\nu)\cdot c(e_a)\cdot
\nabla^{\partial \Omega}_{e_a}\psi$ for an orthonormal basis $e_1, e_2$ for the tangent bundle of $\Sigma$.

 Let $\tilde{M}_{ r}\subset \tilde{M}$ be the region with $\partial
\tilde{M}_{ r}=\Sigma_r$. On $\tilde{M}_r\backslash \Omega$, we
have
\[\begin{split}\int_{\tilde{M}_r\backslash\Omega} (|\nabla\psi|^2-|D\psi|^2)
&=\int_{\partial \Omega} \langle \psi, {D}^{\partial \Omega} \psi+
\frac{1}{2}(k-\langle X,
\nu\rangle)
\psi\rangle\\
&+\int_{\Sigma_r} \langle \nabla_{\nu_r} \psi +c(\nu_r) \cdot
D\psi, \psi \rangle.
\end{split}\]
Adding these up, we obtain
\begin{equation}\label{lich}\begin{split}&\int_{\tilde{M}_r}
|\nabla\psi|^2+\frac{1}{4}\int_\Omega (R+2div
X)|\psi|^2+\frac{1}{2}\int_\Omega X(|\psi|^2)
\\&=\int_{\tilde{M}_r}|D\psi|^2+\int_{\Sigma_r}\langle
\nabla_{\nu}\psi +c(\nu) \cdot D\psi, \psi \rangle.
\end{split}\end{equation}

We claim the left hand side of (\ref{lich}) is always greater than or equal to
\[ \frac{1}{2} \int_{\tilde{M}_r}
|\nabla\psi|^2.\]

This follows from the inequality:
\[|\nabla\psi|^2+|X|^2|\psi|^2+
X(|\psi|^2)\geq 0\] as
\[X(|\psi|^2)=\langle \nabla_X\psi, \psi\rangle
+\langle \psi, \nabla_X \psi\rangle\geq
-2|\nabla_X\psi||\psi|.\]

Thus if we can solve the harmonic spinor equation $D\psi=0$ we obtain
\[\lim_{r\rightarrow\infty} \int_{\Sigma_r}\langle
\nabla_{\nu}\psi +c(\nu) \cdot D\psi, \psi \rangle\geq 0\] and it is known that the limit expression for a constant spinor gives the total mass of $\tilde{M}$.

(\ref{lich}) also implies the following coercive estimates for spinors of compact support
\[\int_{\tilde{M}_r} |D\psi|^2\geq
\frac{1}{2}\int_{\tilde{M}_r} |\nabla\psi|^2\] which is enough to establish the existence of harmonic spinors that are asymptotic to constant spinors at infinity.
\end{proof}

\begin{thm}\label{pos_with_frame}
Given an embedding $i:\Sigma\hookrightarrow N$ into a space-time that satisfies the dominant energy condition. Suppose $\tau$ is admissible, then we have
\[\int_\Sigma h(\Sigma, i_0, \tau, \breve{e}_3) dv_\Sigma \geq \int_\Sigma h(\Sigma, i, \tau, e_3') dv_\Sigma\] where $i_0:\Sigma \hookrightarrow \R^{3,1}$ is the isometric embedding into the Minkowski space given by Theorem B.
\end{thm}

\begin{proof}
Since $\tau$ is  admissible, by (2) and (3) of Definition \ref{adm}, $\Sigma$ bounds a space-like hypersurface $\Omega$ such that Jang's equation over $\Omega$ with boundary value $\tau$ on $\Sigma$ is solvable and the generalized mean curvature $h(\Sigma, i, \tau, e'_3)$ is positive. It follows from Theorem \ref{jang_boundary} that \[\tilde{k}-\langle\widetilde{\nabla}_{\tilde{e}_4} \tilde{e}_4, \tilde{e}_3\rangle+P(\tilde{e}_4, \tilde{e}_3)>0\] on $\widetilde{\Sigma}$, the graph of $\tau$ over $\Sigma$. Take $X$ to be the vector field on $\widetilde{\Omega}$ dual to $\langle \widetilde{\nabla}_{\tilde{e}_4} \tilde{e}_4, \cdot\rangle-P(\tilde{e}_4, \cdot)$, we see that $\tilde{\Omega}$ satisfies the assumption of Theorem \ref{jang_sol} by Remark \ref{shitam}. We can take the projection of $i_0$ onto the standard $\R^3$ slice determined by $t=0$ and denote the image surface by $\widehat{\Sigma}$. The induced metric on $\widehat{\Sigma}$ is then isometric to the metric on the boundary $\widetilde{\Sigma}$ of $\widetilde{\Omega}$. Therefore, by Theorem \ref{jang_sol} we have

\begin{equation}\label{pos}\int_{\widehat{\Sigma}} \hat{k} dv_{\widehat{\Sigma}} \geq \int_{\tilde{\Sigma}} \tilde{k}-\langle \tilde{X}, \tilde{e}_3\rangle dv_{\widetilde{\Sigma}}.\end{equation}  The theorem follows from equation (\ref{int_hat_k}) and  equation (\ref{int_tilde_k}).
\end{proof}

We recall the statement of Theorem A and give the proof:

\bigskip
\noindent\textbf{Theorem A}\,\,\,\textit{Let $N$ be a space-time that satisfies the dominant energy condition. Suppose $i:\Sigma \hookrightarrow N$ is a closed embedded space-like two-surface in $N$ with space-like mean curvature vector $H$. Let $i_0:\Sigma \hookrightarrow \R^{3,1}$ be an isometric embedding into the Minkowski space and let $\tau$ denote the restriction of the time function $t$ on $i_0(\Sigma)$. Let $\bar{e}_4$ be the future-directed time-like unit normal along $i(\Sigma)$ such that
\[\langle H, \bar{e}_4\rangle=\frac{-\Delta \tau}{\sqrt{1+|\nabla \tau|^2}}\] and $\bar{e}_3$ be the space-like unit normal along $\Sigma$ with $\langle \bar{e}_3, \bar{e}_4\rangle=0$ and $\langle H, \bar{e}_3\rangle <0$. Let $\widehat{\Sigma}$ be the projection of $i_0(\Sigma)$ onto $\R^3=\{t=0\}\subset \R^{3,1}$ and $\hat{k}$ be the mean curvature of $\widehat{\Sigma}$ in $\R^3$. If $\tau$ is admissible (see Definition \ref{adm}), then
\[\int_{\widehat{\Sigma}} \hat{k}-\int_\Sigma -\sqrt{1+|\nabla \tau|^2}\langle H, \bar{e}_3\rangle-\alpha_{\bar{e}_3} (\nabla \tau)\] is non-negative.}

\bigskip
\begin{proof}
Because $\tau$ is admissible and the $i_0$ is the unique isometric embedding into $\R^{3,1}$ associated with $i_0$. By equation (\ref{pos}) and equation
(\ref{int_tilde_k}), we have

\begin{equation}\label{pos1}\int_{\widehat{\Sigma}} \hat{k} dv_{\widehat{\Sigma}} \geq \int_\Sigma h(\Sigma, i, \tau, e_3') dv_\Sigma.\end{equation}

By Proposition \ref{can_gauge}, $\int_\Sigma h(\Sigma, i, \tau, e_3') dv_\Sigma\geq \mathfrak{H} (\Sigma, i, \tau)$. Theorem A now follows from Definition \ref{mathfrak_H}.
\end{proof}

Rewriting the integrals, we obtain:
\begin{cor}\label{pos_with_rho}
Given an embedding $i:\Sigma\hookrightarrow N$ into a space-time that satisfies the dominant energy condition. Suppose the mean curvature vector of $\Sigma$ in $N$ is space-like and $\tau$ is admissible, then
\[\mathfrak{H} (\Sigma, i_0, \tau) \geq \mathfrak{H} (\Sigma, i, \tau)\] where $i_0:\Sigma \hookrightarrow \R^{3,1}$ is the isometric embedding into the Minkowski space given by Theorem B.
\end{cor}

\begin{proof} We have from  equation (\ref{int_hat_k}) and equation (\ref{int_hat_k2})
\[\int_{\widehat{\Sigma}} \hat{k} dv_{\widehat{\Sigma}}=\int_\Sigma h(\Sigma, i_0, \tau, \breve{e}_3) dv_\Sigma= \mathfrak{H} (\Sigma, i_0, \tau)\] and

\[\int_\Sigma h(\Sigma, i, \tau, e_3') dv_\Sigma \geq \mathfrak{H} (\Sigma, i, \tau)\] from Proposition \ref{can_gauge}.
\end{proof}

\begin{cor} Under the assumption of Theorem \ref{pos_with_rho}, if the set of admissible $\tau$ is non-empty, then  the quasi-local mass
is non-negative.  It is zero if the embedding $i:\Sigma \hookrightarrow N$ is isometric to $\R^{3,1}$ along $\Sigma$.
\end{cor}

\begin{proof}
The first part follows from the previous corollary. If $i:\Sigma \hookrightarrow N$ is isometric to $\R^{3,1}$ along $\Sigma$, we can take the isometric
embedding $i_0:\Sigma \hookrightarrow \R^{3,1}$, the restriction of the time function $\tau$ will be admissible and all the inequalities become equalities by the uniqueness of $\bar{e}_4$.
\end{proof}

By the boundary gradient estimate of Jang's equation, a constant function is admissible if $\Sigma$ has positive Gauss curvature and the mean curvature vector of $\Sigma$ in $N$ is space-like.

\begin{cor} Under the assumption of Theorem \ref{pos_with_rho}, and supposing $\Sigma$ has positive Gauss curvature, then the quasi-local mass
is non-negative.  It is zero if the embedding $i:\Sigma \hookrightarrow N$ is isometric to $\R^{3,1}$ along $\Sigma$.
\end{cor}

Suppose the minimum is achieved at some $\tau$, we can consider the isometric embedding determined by $\tau$ and define a quasi-local energy momentum vector. This is particularly useful when we have a family of surface $\Sigma_s \hookrightarrow N$, we find the optimal isometric embedding into $\R^{3,1}$ and apply the procedure to get a family of future-directed time-like vectors in $\R^{3,1}$.

\section{The equation of the optimal isometric embedding}

\subsection{Variation of total mean curvature}

Let $\Sigma$ be an orientable closed embedded hypersurface in $\R^{n+1}$. Denote the outward normal by $\nu$ and the mean curvature with respect to $\nu$ by $H$.  We study how the total mean curvature
$\int_{\Sigma} H dv_{\Sigma}$ change with respect to the induced metric $\sigma_{ij}$.

We fix a local coordinate system $u^i$ on $\Sigma$. The variational field $\delta X=Y$ can be decomposed into the tangential and normal part
\[Y=a^k \frac{\partial X}{\partial u^k}+b \nu.\]

We compute the variation of the induced metric $\sigma_{ij}=\langle \frac{\partial X}{\partial u^i}, \frac{\partial X}{\partial u^j}\rangle$.

\begin{equation}\label{del_g}\begin{split}\delta \sigma_{ij}&=\langle \frac{\partial}{\partial u^i}(a^k\frac{\partial X}{\partial u^k}+b\nu), \frac{\partial X}{\partial u^j}\rangle+\langle \frac{\partial X}{\partial u^i}, \frac{\partial}{\partial u^j}(a^k\frac{\partial X}{\partial u^k}+b\nu)\rangle\\
&=\frac{\partial a^k}{\partial u^i}\sigma_{kj}+a^k \Upsilon_{ikj}+b h_{ij} +\frac{\partial a^k}{\partial u^j}\sigma_{ik}+a^k \Upsilon_{jki}+b h_{ij}.\\\end{split}\end{equation} where $\Upsilon_{ikj}=\langle \frac{\partial^2 X}{\partial u^i\partial u^k}, \frac{\partial X}{\partial u^j}\rangle$ is the Christoffel symbol of $\sigma_{ij}$ and $h_{ij}=\langle \frac{\partial \nu}{\partial u^i}, \frac{\partial X}{\partial u^j}\rangle$ is the second fundamental form. Denote by the $\nabla_i$ the covariant derivative with respect to $\frac{\partial}{\partial u^i}$. We solve for
\begin{equation}\label{bh_ij}bh_{ij}=\frac{1}{2}(\delta{\sigma_{ij}}-\nabla_i a^k \sigma_{kj}-\nabla_j a^k \sigma_{ik}).\end{equation}

Next we compute the variation of the mean curvature $H=-\sigma^{ij} \langle \frac{\partial^2 X}{\partial u^i \partial u^j}, \nu\rangle$.
\begin{equation}\label{del_H1}\delta H=h_{ij} \delta \sigma^{ij} -\sigma^{ij}\langle \frac{\partial^2 Y}{\partial u^i \partial u^j}, \nu\rangle-\sigma^{ij}\langle \frac{\partial^2 X}{\partial u^i \partial u^j},  \delta \nu\rangle\end{equation}

We derive
\[\delta \sigma^{ij}=-\sigma^{ik}\delta \sigma_{kl} \sigma^{lj}, \] and
\[\delta \nu=(a^kh_{lk}-\frac{\partial b}{\partial u^l})\sigma^{lj}\frac{\partial X}{\partial u^j}.\]

On the other hand,
\[\begin{split}\langle \frac{\partial^2 Y}{\partial u^i \partial u^j}, \nu\rangle &=\langle \frac{\partial^2 }{\partial u^i \partial u^j}(a^k\frac{\partial X}{\partial u^k}+b\nu), \nu\rangle\\
&=\langle \frac{\partial }{\partial u^i}( \frac{\partial a^k}{\partial u^j}\frac{\partial X}{\partial u^k}+ a^k\frac{\partial^2 X}{\partial u^j \partial u^k}+\frac{\partial b}{\partial u^j}\nu +b\frac{\partial \nu}{\partial u^j}), \nu\rangle\\\end{split}\]

Substitute in
\[\frac{\partial^2 X}{\partial u^j \partial u^k}=\Upsilon_{jk}^l \frac{\partial X}{\partial u^l}-h_{jk}\nu \,\,\text{and}\,\, \frac{\partial \nu}{\partial u^j}=h_{jl} \sigma^{lk}\frac{\partial X}{\partial u^k},\]
and we obtain
\[\begin{split}\langle \frac{\partial^2 Y}{\partial u^i \partial u^j}, \nu\rangle
&=\langle \frac{\partial }{\partial u^i}\left[( \nabla_j a^k +bh_{jl}\sigma^{lk})\frac{\partial X}{\partial u^k}+(\frac{\partial b}{\partial u^j}-a^kh_{jk} )\nu\right], \nu\rangle\\
&=-h_{ik}( \nabla_j a^k+bh_{jl}\sigma^{lk})+\frac{\partial}{\partial u^i}(\frac{\partial b}{\partial u^j}-a^lh_{jl} )\end{split}.\]

Plug these into (\ref{del_H1}) and we arrive at

\[\begin{split}\delta H&=-\sigma^{ik}\sigma^{lj}h_{ij} \delta \sigma_{kl}-\Delta b+\sigma^{ij} h_{ik} \nabla_j a^k +bg^{ij} \sigma^{lk} h_{ik} h_{jl}+ \sigma^{ij} \nabla_i (a^k h_{jk}).\end{split}\]

We plug (\ref{bh_ij}) into this equation and obtain
\[\delta H=-\frac{1}{2}\sigma^{ik}\sigma^{lj}h_{ij} \delta \sigma_{kl}-\Delta b+\sigma^{ij}\nabla_i (a^k h_{jk}).\]

\begin{pro} Let $\Sigma$ be a closed embedded hypersurface in $\R^{n+1}$.
The variation of the total mean curvature with respect to a metric deformation is

\[\delta \int_{\Sigma} H dv_{\Sigma}=\frac{1}{2}\int_\Sigma (H \sigma^{ij}-\sigma^{ik} \sigma^{jl} h_{kl})\delta{ \sigma_{ij}} dv_{\Sigma}.\]
\end{pro}
\begin{cor}
If $X_s:\Sigma\rightarrow \R^{n+1}$ is a smooth family of isometric embedding of a compact $n$-manifold $\Sigma$, then the total mean curvature is a constant.
\end{cor}

\subsection{The variational equation}
In this section, $\sigma_{ab}$ denotes the metric on a two-surface $\Sigma$ which satisfies the assumption in Theorem A. Recall the metric on the projection $\widehat{\Sigma}$ is $\hat{\sigma}_{ab}=\sigma_{ab}+\tau_a\tau_b$. The metric $\sigma_{ab}$ is fixed for the isometric embedding and thus $\delta\hat{\sigma}_{ab}=\delta(\tau_a\tau_b)$.

The quasi-local mass expression we try to minimize is

\[\begin{split}&\Xi=\int_{\widehat{\Sigma}} \hat{k} dv_{\widehat{\Sigma}}-\int_\Sigma \left[\sqrt{1+|\nabla\tau|^2}\cosh\theta|{H}|-\nabla\tau\cdot \nabla \theta -V \cdot \nabla \tau \right]dv_\Sigma.\end{split}\]
where
$\sinh \theta =\frac{-\Delta \tau}{|{H}|\sqrt{1+|\nabla \tau|^2}}$, and $V$ is the tangent vector on $\Sigma$ that is dual to the connection one-form $\alpha_{\hat{e}_3}$ determined by  $\hat{e}_3=-\frac{H}{|H|}$. This is an expression that determined by $\sigma_{ab}$ and the mean curvature vector $H$.

We can take $X=( x(u^a), y(u^a), z(u^a), \tau (u^a)):\Sigma \hookrightarrow \R^{3,1}$ and $\widehat{X}=( x(u^a), y(u^a), z(u^a)): \Sigma \hookrightarrow \R^3$.  Notice that $dv_{\widehat{\Sigma}}=\sqrt{\det\hat{\sigma}_{ab}} du^1 du^2$ and $dv_{{\Sigma}}=\sqrt{\det{\sigma}_{ab}} du^1 du^2$.
Recall from the last section the variation of $\int_{\widehat{\Sigma} }\hat{k} dv_{\widehat{\Sigma}}$ is

\[\delta \int_{\widehat{\Sigma}} \hat{k} dv_{\widehat{\Sigma}}=\int_{\widehat{\Sigma}}(\widehat{H} \hat{\sigma}^{ab} -\hat{\sigma}^{ac} \hat{\sigma}^{bd} \hat{h}_{cd})\tau_a (\delta\tau)_b dv_{\widehat{\Sigma}}.\]

Integration by parts and recall that the tensor $\widehat{H} \hat{\sigma}^{ab} -\hat{\sigma}^{ac} \hat{\sigma}^{bd} \hat{h}_{cd}$ is divergence free
on $\hat{\Sigma}$, we obtain

\begin{equation}\begin{split}\delta \int_{\widehat{\Sigma}} \hat{k} dv_{\widehat{\Sigma}}
&=-\int_{\widehat{\Sigma}}(\widehat{H} \hat{\sigma}^{ab} -\hat{\sigma}^{ac} \hat{\sigma}^{bd} \hat{h}_{cd})\widehat{\nabla}_b\widehat{\nabla}_a \tau {\delta\tau} dv_{\widehat{\Sigma}}\end{split}\end{equation}

The relation between the Hessians of $\tau $ on $\Sigma$ and $\widehat{\Sigma}$ is given by
\[\widehat{\nabla}_b\widehat{\nabla}_a\tau=\frac{1}{1+|\nabla\tau|^2} \nabla_b\nabla_a \tau.\]

Since we also have

\[\frac{\det \hat{\sigma}}{\det \sigma}=1+|\nabla \tau|^2\]

\begin{equation}\begin{split}\delta \int_{\widehat{\Sigma}} \hat{k} dv_{\widehat{\Sigma}}
&=-\int_{\widehat{\Sigma}}(\widehat{H} \hat{\sigma}^{ab} -\hat{\sigma}^{ac} \hat{\sigma}^{bd} \hat{h}_{cd})\frac{\nabla_b\nabla_a \tau}{\sqrt{1+|\nabla\tau|^2}}{\delta\tau} dv_{{\Sigma}}\end{split}\end{equation}

Now \[ \begin{split}&\delta\int_\Sigma \left[\sqrt{1+|\nabla\tau|^2}\cosh\theta|{H}|-\nabla\tau\cdot \nabla \theta \right]dv_\Sigma\\
&=\int_\Sigma \left[(1+|\nabla\tau|^2)^{-1/2}\nabla \tau \cdot \nabla \delta \tau \cosh\theta|{H}|+\sqrt{1+|\nabla \tau|^2}\sinh\theta \delta \theta|{H}| \right]dv_\Sigma\\
&-\int_\Sigma (\nabla\delta \tau \cdot \nabla\theta+\nabla \tau\cdot \nabla\delta \theta ) dv_\Sigma.\end{split}\]

Substitute $\sinh \theta =\frac{-\Delta \tau}{|{H}|\sqrt{1+|\nabla \tau|^2}}$ and integrate by parts, we obtain

\[ \begin{split}&\delta\int_\Sigma \left[\sqrt{1+|\nabla\tau|^2}\cosh\theta|{H}|-\nabla\tau\cdot \nabla \theta \right]dv_\Sigma\\
&=\int_\Sigma (\frac{\nabla\tau}{\sqrt{1+|\nabla\tau|^2}} \cosh\theta|{H}|-\nabla\theta)\cdot \nabla \delta\tau dv_\Sigma\end{split}\]

\begin{pro}
The variation of $\Xi$ with respect to $\tau$ is
\[-\int_{{\Sigma}}(\widehat{H} \hat{\sigma}^{ab} -\hat{\sigma}^{ac} \hat{\sigma}^{bd} \hat{h}_{cd})\frac{\nabla_b\nabla_a \tau}{\sqrt{1+|\nabla\tau|^2}}{\delta\tau} dv_{{\Sigma}}+\int_\Sigma div_\Sigma (\frac{\nabla\tau}{\sqrt{1+|\nabla\tau|^2}} \cosh\theta|{H}|-\nabla\theta-V)\cdot  \delta\tau dv_\Sigma \]
\end{pro}

Therefore, the equation for the minimizing isometric embedding is
\begin{equation}-(\widehat{H}\hat{\sigma}^{ab} -\hat{\sigma}^{ac} \hat{\sigma}^{bd} \hat{h}_{cd})\frac{\nabla_b\nabla_a \tau}{\sqrt{1+|\nabla\tau|^2}}+ div_\Sigma (\frac{\nabla\tau}{\sqrt{1+|\nabla\tau|^2}} \cosh\theta|{H}|-\nabla\theta-V)=0\end{equation} with $\sinh \theta =\frac{-\Delta \tau}{|{H}|\sqrt{1+|\nabla \tau|^2}}$.

\end{document}